\documentclass[10pt, a4paper, notitlepage]{article}

\usepackage{amssymb}
\usepackage{amsmath}
\usepackage{amsthm}
\usepackage{natbib}
\usepackage[small]{titlesec}
\usepackage{color}

\newcommand{\diag}{{\text{diag}}}

\newcommand{\PP}{\mathbb{P}}

\newcommand{\lle}{\,\,{\lesssim}\,\,}

\newcommand{\eps}{\varepsilon}

\newcommand{\FF}{\mathcal{F}}

\newcommand{\RR}{\mathbb{R}}

\newcommand{\NN}{\mathbb{N}}

\newcommand{\HHH}{\mathbb{H}}

\renewcommand{\phi}{\varphi}

\newcommand{\given}{\,|\,}

\newtheorem{thm}{Theorem}[section]
\newtheorem{lem}[thm]{Lemma}

\theoremstyle{definition}

\newenvironment{prf}{\vspace{1ex}\begin{proof}[\bf Proof]}{\end{proof}\vspace{2ex}}

\begin{document}

\title{Full adaptation to smoothness using randomly
truncated series priors with Gaussian coefficients and inverse gamma scaling\\[.5cm]}

\author{
Jan van Waaij and Harry van Zanten \footnote{Research  funded by the Netherlands Organization for Scientific Research (NWO)}
\\[2ex]
 Korteweg-de Vries Insitute for Mathematics\\
 University of Amsterdam \\[1ex]
 \texttt{j.vanwaaij@uva.nl, hvzanten@uva.nl}
}

\date{December 5, 2016 \ (revised version)}

\maketitle

\begin{abstract}
We study random series priors for estimating a functional parameter $f \in L^2[0,1]$.  
We show that with a series prior with random truncation, Gaussian coefficients, and 
inverse gamma multiplicative scaling, it is possible to achieve posterior 
contraction at optimal rates and adaptation to arbitrary degrees of smoothness.  
 We present general results that can be combined with existing rate of contraction 
 results for various nonparametric estimation problems. We give concrete examples 
 for signal estimation in white noise and drift estimation for a one-dimensional SDE.
\end{abstract}

\numberwithin{equation}{section}

%\begin{center}
%{\sc {\bf Running head:} Spectral theory  for the fBm}
%\end{center}

\section{Introduction}

In Bayesian function estimation, a common approach to putting a prior 
distribution on a function $f$ of interest, for instance a regression function
in nonparametric regression models or a drift function in diffusion models,  is 
to expand the function in a particular basis
and to endow the coefficients in the expansion with prior weights.
For computational or other reasons the series is often truncated after finitely many terms, 
and the truncation level is endowed with a prior as well. The coefficients
in the expansion are often chosen to be independent under the prior and 
distributed according to some given probability density. 

It is of interest to understand whether, in addition to their attractive
conceptual and computational aspects, nonparametric priors of this 
type enjoy favourable theoretical properties as well. 
Examples of 
papers in which this  was studied for various families of series priors include \cite{zhao2000}, 
\cite{shenwasserman2001}, \cite{dejonge2012},
\cite{rivoirard2012}, \cite{arbel2013}, \cite{shen2015}.
The results in these papers show that when appropriately constructed, 
random series priors can yield posteriors that contract at optimal rates 
and that adapt automatically to the smoothness of the function that is being estimated.

To ensure that the nonparametric Bayes procedure not only adapts to smoothness, but  is also flexible 
with respect to the multiplicative scale 
of the function of interest, a multiplicative hyperparameter with an independent
prior distribution is often employed as well. Theoretically this is usually not 
needed for an optimal concentration rate of the posterior, but it  can greatly improve performance in practice. {See for instance 
  \cite{meulen2014}, where it is explained why it is } 
%A 
computationally attractive in certain settings % choice is 
 to use Gaussian priors on the 
series coefficients in combination with a multiplicative (squared) scaling parameter 
with an inverse gamma prior. For a given truncation level, the prior is conjugate and 
allows for  posterior computations using standard Gibbs sampling. 
The existing theoretical results do not cover this important case however. This is mainly due 
to the fact that essentially, the available rate of contraction 
theorems for series priors require that hyper priors have (sub-) exponential tails, 
which excludes the inverse gamma distribution. { (For example the second part of condition (A2) of \cite{shen2015} is not satisfied in our setting.)} 
The theoretical properties of random series priors with inverse gamma scaling have therefore
remained unexplored. With this paper we intend to fill this gap.

Concretely, we consider statistical models in which the unknown object of 
interest is a square integrable function $f$ on $[0,1]$. We endow this function with a prior 
that is hierarchically specified as follows:
\begin{equation}\label{eq: prior}
\begin{split}
J & \sim \text{Poisson or geometric},\\
s^2 & \sim \text{inverse gamma},\\
f \given s, J & = \sum_{j \le J} f_j\psi_j, \quad \text{with} \  
(f_1, \ldots, f_J) \sim N(0, \diag(s^2j^{-1-2\alpha})_{j \le J}),
\end{split}
\end{equation}
where $(\psi_j)$ is a fixed orthonormal basis of $L^2[0,1]$ and $\alpha > 0$ is a hyperparameter.
(In fact, we will consider a somewhat broader class of hyper priors on $J$ and $s^2$, 
see Section \ref{sec: prior}.) %{\color{cyan} Every element \(f\) of \(L^2[0,1]\) can be represented as a Fourier series in this basis and the decay of the coefficients determine the smoothness of the function and how well \(f\) is approximated when it is projected on the linear span of the first \(j\) basis elements.} 
%Simulation studies in   \cite{meulen2014}, who considered nonparametric estimation of the 
%drift in  a stochastic differential equation (SDE) model,  show that such a prior is  indeed 
%quite flexible,  with respect both to the scale and to the smoothness of the unknown function. 
In this paper we prove that this prior enjoys very favourable theoretical properties as well.
We derive optimal posterior contraction rates and adaptation up
%In particular, we prove that with this prior it is possible to achieve adaptation 
to arbitrarily high degrees of smoothness.

In recent years,  general rate of contraction theorems have been derived for 
a variety of nonparametric statistical problems. Roughly speaking, such theorems give sufficient conditions
for having a certain rate of contraction in terms of (i) the amount of mass that the prior 
gives to neighbourhoods of the true function and (ii) the existence of  growing subsets of 
the support of the prior, so-called sieves, that contain all but an exponentially small amount of the prior mass 
and whose metric entropy is sufficiently small. The statements of our main theorem match the conditions
of these existing general results. This means that we automatically obtain results for different statistical settings,
including for instance  signal estimation in white noise and drift estimation for SDEs.

A simple but important observation that we make in this paper is that to obtain sharp rates 
for the priors we consider, it is necessary to use  versions of the general contraction rate  theorems 
that give entropy conditions on the intersection of the sieves with balls around the true function,
as can be found  for instance in \cite{GGV}, \cite{meulen2006} and \cite{GVnoniid}. 
As remarked in these papers, it is in many nonparametric problems sufficient to consider 
only the entropy of the sieves themselves, without intersecting them with a ball around the truth. For the priors
we  consider in this paper however, which in some sense are  finite-dimensional in nature in certain regimes, 
this is not the case. It turns out that since the inverse gamma 
has polynomial tails, we need to make the sieves relatively large in order to ensure that they receive
sufficient prior mass. Without intersecting them with a small ball around the truth, this would 
make their entropy too large, or even infinite. 

The proof of our main results indicate 
that the good adaptation properties of series priors like \eqref{eq: prior}
are really due to the fact that {\em both} the truncation level $J$ {\em and} the scaling constant $s$
are random. If the true function that is being estimated is relatively smooth, the prior
can approximate it well by letting $J$ be small. If it is relatively rough however, the prior can adapt 
to it by letting $J$ be essentially infinite, or very large, to pick up all the fluctuations. The 
correct bias-variance trade-off is in that case achieved automatically by adapting the multiplicative scale.
In some sense, priors like \eqref{eq: prior} can switch with sufficient probability between 
being essentially finite-dimensional, and being essentially infinite-dimensional. In combination 
with a random multiplicative scale, this gives them the ability to adapt to all levels of smoothness. %{\color{cyan} When \(J\) is infinite or very large, \(\alpha\) determines the baseline smoothness of the prior. However, we obtain adaptivity for every \(\beta>0\) regardless of the choice of \(\alpha\).}

The remainder of the paper is organized as follows. In the next section we describe 
in detail the class of priors we consider. %, which is a bit larger than just \eqref{eq: prior}, 
%in the sense that we allow larger classes of hyper priors on $J$ and $s$. 
In Section \ref{sec: main} we present the main results of the paper, which give 
bounds on the amount of mass that the priors give to $L^2$-neighbourhoods of functions with a given 
degree of (Sobolev-type) smoothness, and the existence of appropriate sieves within the support 
of the prior. In Section \ref{sec: app} we link these general  theorems to existing rate of contraction results
for two different SDE models, to  obtain concrete contraction results for signal estimation in white
noise and drift estimation of a one-dimensional SDE with priors of the form \eqref{eq: prior}. 
The proofs of the main results are given in Sections \ref{sec: proof1} and \ref{sec: proof2}.

%\subsection{Notation}
%
%$\lle$, $\gtrsim$, $\asymp$, $\|\cdot\|_2$
%
%$L^2[0,1]$, $H^\beta[0,1]$, $B(f_0, \eps)$. 
%
%
%$N(\eps, \FF, \|\cdot\|)$.

\section{Prior model}
\label{sec: prior}

We consider problems in which the unknown function of interest (e.g.\ a drift function of an SDE, a
signal observed in noise, \ldots) 
is a square integrable function on $[0,1]$, 
i.e.\ an element of  $L^2[0,1] = \{f: [0,1] \to \RR: \|f\|_2 < \infty\}$, 
where the $L^2$-norm is as usual defined by 
$\|f\|^2_2 = \int_0^1 f^2(x)\,dx$. %To construct a prior on this space 
 We fix an arbitrary orthonormal basis $(\psi_j)$ of $L^2[0,1]$ (for instance the standard Fourier basis). { Every element of \(f\in L^2[0,1]\) can be represented as a series \(f=\sum_j \langle f,\psi_j\rangle \psi_j\) where the convergence is in the \(L^2\)-norm and by the Plancherel formula \(\|f\|_2^2=\sum_{j}|\langle f,\psi_j\rangle|^2\). Finite series \(\sum_{j\le J}\langle f,\psi_j\rangle \psi_j\) approximate \(f\) and the quality of this approximation depends on the decay of the coefficients \(\langle f,\psi_j\rangle\), which also determines the ``smoothness'' of the function. The class of \(\beta\)-Sobolev smooth functions \(H^\beta[0,1]\) is given by all \(f\in L^2[0,1]\) for which the \(\beta\)-Sobolev norm \(\|f\|_\beta:=\sqrt{\sum_{j}k^{2\beta}|\langle f,\psi_j\rangle|^2}\) 
is finite. If $\psi_j$ is the classical Fourier series basis, 
these are the classical\(\beta\)-Sobolev spaces.}

% We use the basis $(\psi_j)$ of $L^2[0,1]$ to construct a prior on \(L^2[0,1]\). In order to do so, fix a number $\alpha > 0$.}
We define a series prior on a function $f \in L^2[0,1]$ through a hierarchical 
scheme which involves a prior on the point $J$ at which the series is truncated, 
a prior on the multiplicative scaling constant $s$ and conditionally on 
$s$ and $J$,  a series prior with Gaussian coefficients on $f$. 

Specifically, the prior on $J$ is defined through a probability mass
function $p$ that is assumed to satisfy, for constants $C, C' >0$, 
\begin{equation}\label{eq: p}
p(j)  \gtrsim e^{-Cj\log j}, \qquad \sum_{i > j} p(i) \lle e^{-C'j}
\end{equation}
for all $j \in \NN$. (As usual, $a \lle b$ or $b \gtrsim a$ means that $a \le cb$ for 
an irrelevant constant $c >0$.)  This includes for instance the cases of a Poisson or a geometric prior on  $J$. 
For the scaling parameter  we assume that the density  $g$ 
of $s^2$  is positive and continuous and satisfies, for some $q < -1$ and $C'' >0$, 
\begin{equation}\label{eq: g}
g(x) \gtrsim e^{-C''/x} \quad \text{near $0$}, \qquad
g(x) \gtrsim x^{q} \quad \text{near $\infty$}.
\end{equation}
Hence in particular, the popular and 
computationally convenient choice of an inverse gamma prior on $s^2$ is included in our setup.
The full prior {\(\Pi\) is then specified as follows:}
\begin{align}
\label{eq: p1} J & \sim p\\
\label{eq: p2} s^2 & \sim g\\
\label{eq: p3} f \given s, J & \sim s \sum_{j=1}^J j^{-1/2-\alpha} Z_j \psi_j,
\end{align}
where {\(\alpha\) is a positive constant which determines the baseline smoothness of the prior,} $p$ satisfies \eqref{eq: p}, $g$ satisfies \eqref{eq: g} and the $Z_j$
are independent standard Gaussians. 

%The positive hyperparameter $\alpha$ can be chosen freely. Very roughly speaking, it 
%determines the ``baseline regularity'' of the prior, in the sense that it determines 
%the rate of decay of the Fourier coefficients of $f$ relative to the basis $(\psi_j)$. 
%However, this intuition only makes sense if $J$ is infinite or very large. 
%The whole point of the construction of the prior is in fact that it is flexible and that with 
%sufficient  probability, the parameters $J$ and $s$ can take values that change the 
%character of the prior. Indeed, if $J$ is small, it effectively sets the higher Fourier
%coefficients to $0$, making $f$ ``smoother''. Similarly, if $s$ is smaller, it also has 
%the effect of ``smoothing'' $f$ in a sense. Our main result below confirms that we  indeed
%have the desired flexibility, in the sense  that we can estimate functions of arbitrary smoothness
%at optimal rates with this prior, irrespective of our initial choice of the hyperparameter $\alpha$. 

\section{Main results}
\label{sec: main}

Our main abstract result gives properties of the truncated 
series prior that link directly to the conditions of existing 
general theorems for posterior contraction in a variety of statistical settings.
Combined with such existing results, we  obtain concrete 
results for, for instance,  signal estimation in white noise,
drift estimation in diffusion models, et cetera. We give concrete 
examples in the next section.

%For $\beta > 0$, define 
%\[
%H^\beta[0,1] = \Big\{f = \sum f_j\psi_j : \sum f^2_j j^{2\beta} < \infty\Big\}.
%\]
%Note that if $(\psi_j)$ is the ordinary Fourier basis, this is the usual $L^2$-Sobolev 
%space of regularity $\beta$.
 As usual, if $\FF$ is a { subset of a normed vector} space  %that is endowed 
 with norm $\|\cdot\|$, then we denote by $N(\eps, \FF, \|\cdot\|)$
the minimal number of balls of $\|\cdot\|$-radius $\eps$ needed to cover the set $\FF$. 

%\bigskip

\begin{thm}\label{thm: main}
Let the prior $\Pi$ on $f$ be as defined in \eqref{eq: p1}--\eqref{eq: p3},
with $\alpha > 0$ and $p$ and $g$ satisfying \eqref{eq: p}--\eqref{eq: g}. Let $f_0 \in H^\beta[0,1]$ for $\beta > 0$. 
Then there exists a constant $c > 0$
such that for every $K > 1$, there exist $\FF_n \subset L^2[0,1]$ such that with
\begin{equation}\label{eq: epsn}
\eps_n =c \Big(\frac{n}{\log n}\Big)^{-\beta/(1+2\beta)}, 
\end{equation}
we have
\begin{align}
\label{eq: pm} \Pi(f:\|f - f_0\|_2 \le \eps_n) & \ge e^{-n\eps^2_n},\\
\label{eq: rm} \Pi(f \not \in \FF_n)  &\le e^{-Kn\eps_n^2},\\
\label{eq: en} \sup_{\eps > \eps_n}\log N(a\eps,  \{f \in \FF_n : \|f-f_0\|_2 \le \eps\}, \|\cdot\|_2)  & \lle n\eps^2_n,
\end{align}
for all $a \in (0,1)$.
\end{thm}

%\bigskip

The proof of the theorem is given in Section \ref{sec: proof1}. 
The result matches with the sufficient conditions of existing posterior contraction theorems,
provided that the relevant statistical distance-type quantities (e.g.\ Hellinger, Kullback-Leibler, \ldots)
in the model can be appropriately linked to the $L^2$-norm on the parameter $f$. In the next 
section we give two concrete SDE-related examples, which motivated the present study. 

Theorem \ref{thm: main} shows that with truncated series priors of the type \eqref{eq: p1}--\eqref{eq: p3}
we can have adaption to arbitrary degrees of smoothness in certain function estimation problems, and achieve posterior contraction 
rates that are optimal up to a logarithmic factor. 
Inspection of the proof of Theorem \ref{thm: main} shows that in the range $\beta \le \alpha + 1/2$, 
i.e.\ if the ``baseline smoothness'' $\alpha$ of the prior happens to have been chosen large enough 
relative to the smoothness $\beta$ of the true function, then we actually get the optimal rate $n^{-\beta/(1+2\beta)}$
without additional logarithmic factors. This is true under a slightly stronger condition on  the prior on 
the cut-off point $J$. Instead of \eqref{eq: p}, we need to assume that for constants $C, C' > 0$ it holds that  
\begin{equation}\label{eq: pp}
p(j)  \gtrsim e^{-Cj}, \qquad \sum_{i > j} p(i) \lle e^{-C'j}
\end{equation}
for all $j \in \NN$. This means that 
the prior on $J$ can still be geometric, but that the Poisson prior on $J$ is excluded.

%\bigskip

\begin{thm}\label{thm: main2}
Let the prior $\Pi$ on $f$ be as defined in  \eqref{eq: p1}--\eqref{eq: p3},
with $\alpha > 0$ and $p$ and $g$ satisfying \eqref{eq: pp} and \eqref{eq: g}. 
Let $f_0 \in H^\beta[0,1]$ for $0< \beta \le \alpha + 1/2$. Then there exists a constant $c > 0$
such that for every $K > 1$, there exist $\FF_n \subset L^2[0,1]$ such that with
\begin{equation}\label{eq: epsn2}
\eps_n =c n^{-\beta/(1+2\beta)}, 
\end{equation}
we have
\begin{align}
\label{eq: pm2} \Pi(f: \|f - f_0\|_2 \le \eps_n) & \ge e^{-n\eps^2_n},\\
\label{eq: rm2} \Pi(f \not \in \FF_n)  &\le e^{-Kn\eps_n^2},\\
\label{eq: en2} \sup_{\eps > \eps_n}\log N(a\eps,  \{f \in \FF_n : \|f-f_0\|_2 \le \eps\}, \|\cdot\|_2)  & \lle n\eps^2_n,
\end{align}
for all $a \in (0,1)$.
\end{thm}

%\bigskip

The proof of this theorem is given in Section \ref{sec: proof2}.

\section{Specific statistical settings}
\label{sec: app}

\subsection{Detecting a signal in Gaussian white noise}

Suppose we observe a sample path $X^{(n)}=(X^{(n)}_t: t \in [0,1])$ of stochastic process 
satisfying the SDE
\[
dX^{(n)}_t = f_0(t)\,dt + \frac1{\sqrt n}\,dW_t, 
\]
where $W$ is a standard Brownian motion and $f_0 \in L^2[0,1]$ is an unknown signal. 
To make inference about the signal we endow it with the truncated series prior $\Pi$
described in Section \ref{sec: prior} and we compute the corresponding posterior 
$\Pi(\cdot \given X^{(n)})$. Theorem 3.1 of \cite{meulen2006} or Theorem 6 of \cite{GVnoniid}, 
combined by our main result Theorem \ref{thm: main}, imply that if $f_0 \in H^\beta[0,1]$
for $\beta > 0$, then we have the posterior contraction
\[
\Pi(f: \|f-f_0\|_2 >M_n (n/\log n)^{-\beta/(1+2\beta)} \given X^{(n)}) \overset{P_{f_0}}{\to} 0
\]
for all $M_n \to \infty$, where the convergence is in probability under the true 
model corresponding to the signal $f_0$.

\subsection{Estimating the drift of an ergodic diffusion} 

Suppose we observe a sample path $X^{(T)} = (X_t: t \in [0,T])$ of an 
ergodic one-dimensional diffusion satisfying the SDE
\[
dX_t = b_0(X_t)\,dt + \sigma(X_t)\,dW_t, \qquad X_0 = 0, 
\]
where $W$ is a standard Brownian motion, $\sigma: \RR \to \RR$ is a know 
continuous function that is bounded away from $0$, and $b_0: \RR \to \RR$ is a continuous function 
that satisfies the appropriate conditions to guarantee that the  SDE indeed generates an ergodic diffusion
(see for instance \cite{kallenberg2002}). The goal is to estimate the 
restriction $b_0|_{[0,1]}$ of $b_0$ to the interval $[0,1]$.

The likelihood for this model, given by Girsanov's formula (e.g.\ \cite{liptser2001}),  
factorizes into a factor involving only the drift
on the interval $[0,1]$ and a factor involving only the restriction of the drift 
to the complement $\RR\backslash [0,1]$. 
As a result, since we are only interested in the drift on $[0,1]$, we can effectively
assume that it is known outside $[0,1]$ and we only have to put a prior on the restriction of 
the drift to $[0,1]$. We endow this with the truncated series prior $\Pi$
described in Section \ref{sec: prior} and we compute the corresponding posterior 
$\Pi(\cdot \given X^{(T)})$. Theorem 3.3 of \cite{meulen2006}
and  Theorem \ref{thm: main} then imply that if $b_0|_{[0,1]} \in H^\beta[0,1]$
for $\beta > 0$, then we have the posterior contraction
\[
\Pi(b: \|b-b_0\|_2 >M_T (T/\log T)^{-\beta/(1+2\beta)} \given X^{(T)}) \overset{P_{b_0}}{\to} 0
\]
as $T \to \infty$ for all $M_T \to \infty$, where the convergence is in probability under the true 
model corresponding to the drift function $b_0$.

\section{Proof of Theorem \ref{thm: main}}
\label{sec: proof1}

\subsection{Prior mass}

The following theorem implies that \eqref{eq: pm} holds with $\eps_n$ as specified.

\begin{thm}\label{thm: pm1}
Let the prior $\Pi$ on $f$ be defined according to \eqref{eq: p1}--\eqref{eq: p3},  
with $\alpha > 0$ and $p$ and $g$ satisfying \eqref{eq: p}--\eqref{eq: g},
and let $f_0 \in H^\beta[0,1]$ for $\beta > 0$. 
Then, for a constant $C > 0$, it holds that  
\[
-\log \Pi(f: \|f-f_0\|_2 \le 2\eps) \le C \eps^{-1/\beta}\log 1/\eps,
\]
for all $\eps> 0$ small enough.
\end{thm}

\begin{prf}
Recall that $s^2$ has density $g$ under the prior. Hence, by conditioning we see that 
the  probability of interest  is bounded from below by 
\[
p\big(\big\lfloor(\eps/\|f_0\|_\beta)^{-1/\beta}\big\rfloor\big) 
\int_{0}^{\infty} \Pi\Big(\Big\|{\sqrt{\eta}} \sum_{j=1}^{\big\lfloor(\eps/\|f_0\|_\beta)^{-1/\beta}\big\rfloor} j^{-1/2-\alpha} Z_j \psi_j  - f_0\Big\|_2 \le 2\eps\Big)g(\eta)\,d\eta,
\]
%{\color{cyan}with \(\Pi_{J_0,\eta}\) the prior conditioned on \(J=J_0:=\big\lfloor(\eps/\|f_0\|_\beta)^{-1/\beta}\big\rfloor\) en \(s^2=\eta\), which is given by \eqref{eq: p3}. }

Now suppose first that $1+2\alpha-2\beta \le 0$. Then by 
Lemma \ref{lem: sb}, the preceding is further  lower bounded by 
\[
\exp\big(-C_1\eps^{-1/\beta}\log 1/\eps\big) 
p\big(\big\lfloor(\eps/\|f_0\|_\beta)^{-1/\beta}\big\rfloor\big) \int_{\eps^{1/\beta}}^{2\eps^{1/\beta}}g(\eta)\,d\eta 
\]
for a constant $C_1> 0$.  By the assumptions on $p$ and $g$ this is 
bounded from below by a constant times $\exp(-C_2\eps^{-1/\beta}\log 1/\eps)$
for $\eps$ small enough, for some constant \(C_2>0\). 
  
In the other case $1+2\alpha-2\beta > 0$ we restrict the integral over $\eta$ to 
a different region to obtain instead the lower bound 
\[
\exp\big(-C_1\eps^{-1/\beta}\log 1/\eps\big) 
p\big(\big\lfloor(\eps/\|f_0\|_\beta)^{-1/\beta}\big\rfloor\big) 
\int_{\eps^{(2\beta - 2\alpha)/\beta}}^{2\eps^{(2\beta - 2\alpha)/\beta}}g(\eta)\,d\eta 
\]
for some $C_1> 0$. For $\alpha < \beta \le \alpha + 1/2$ 
the assumptions on $p$ and on the behaviour of $g$ near $0$  ensure again that this 
is 
bounded from below by a constant times $\exp(-C_2\eps^{-1/\beta}\log 1/\eps)$
for $\eps$ small enough. For the range $\beta < \alpha$ this holds as well, 
by the 
the assumptions on $p$ and on the behaviour of $g$ near $\infty$. 
{When \(\alpha=\beta\) we use the lower bound 
\[
\exp\big(-C_1\eps^{-1/\beta}\log 1/\eps\big) 
p\big(\big\lfloor(\eps/\|f_0\|_\beta)^{-1/\beta}\big\rfloor\big) 
\int_{1}^{C_3}g(\eta)\,d\eta.\]
Again by the behaviour of \(g\) near infinity, the integral on the right is positive for \(C_3\) big enough and the desired lower bound holds by the assumption on \(p\).} 
\end{prf}

\begin{lem}\label{lem: cb}
Let $Z_1, Z_2, \ldots$ be independent and standard normal.
There exists a universal constant $K> 1$  
such that for every $s > 0$, $\eps > 0$, $J \in \NN$ and $a \in \ell^2$,
\[
-\log \PP\Big(\Big\|s \sum_{j=1}^J a_j Z_j \psi_j \Big\|_2 \le \eps\Big) \le
2 J  \log \Big(K \vee \frac{s\|a\|_2}{\eps}\Big).
\]
\end{lem}

\begin{prf}
Since the $\psi_j$ form an orthonormal basis, the probability we have to lower bound 
equals
\[
\PP\Big(s^2 \sum_{j=1}^J a^2_j Z^2_j \le \eps^2\Big) \ge 
\PP\Big( \max_{j \le J} |Z_j|   \le \frac{\eps}{s\|a\|_2}\Big) = 
\Big(\PP\Big( |Z_1|   \le \frac{\eps}{s\|a\|_2}\Big)\Big)^J.
\]
If ${\eps}/({s\|a\|_2}) \ge \xi_{3/4}$, with $\xi_p$ the $p$-quantile of the 
standard normal distribution, { \(\PP( |Z_1|   \le \eps/(s\|a\|_2))\ge 1/2.\)} 
% the righthand side is at least $(1/2)^J$. 
In the other case, it is at least  { \(\phi(\xi_{3/4}) \times {2\eps}/({s\|a\|_2})\)},
with $\phi$ the standard normal density.
So in either case, it is at least a constant $C \in (0,1)$ times 
$1 \wedge {\eps}/({s\|a\|_2})$. It follows that 
\begin{align*}
\log \PP\Big(\Big\|s \sum_{j=1}^J a_j Z_j \psi_j \Big\|_2 \le \eps\Big) & \ge 
J \log C+ J \log \Big(1 \wedge \frac{\eps}{s\|a\|_2}\Big)\\
& \ge 2 J \log \Big(C \wedge \frac{\eps}{s\|a\|_2}\Big).
\end{align*}
This implies the statement of the lemma.
\end{prf}

\begin{lem}\label{lem: sb}
Let $Z_1, Z_2, \ldots$ be independent and standard normal.
Let \(\beta>0\) and \(f_0\in H^\beta[0,1]\) be given. 
There exists a  constant $K > 1$ such that for all 
$\eps, s, \alpha > 0$ and 
$J \ge (\eps/\|f_0\|_\beta)^{-1/\beta}$, 
\[
-\log \PP\Big(\Big\|s \sum_{j=1}^J j^{-1/2-\alpha} Z_j \psi_j  - f_0\Big\|_2 \le 2\eps\Big) \le
2 J  \log \Big(K \vee \frac{s}{\eps}\Big) + 
\frac{\|f_0\|^2_\beta}{s^2}J^{{(1+2\alpha-2\beta) \vee 0}}.
\]

\end{lem}

\begin{prf} {For  fixed \(J,s\), the sum
%We can view 
$s \sum_{j=1}^J j^{-1/2-\alpha} Z_j \psi_j$ is a centered Gaussian random element in $L^2[0,1]$ and has a
%As such is has an associated 
 reproducing kernel Hilbert space (RKHS),} which is the space $\HHH^{s,J}$ 
of all functions $h = \sum_{j\le J} h_j\psi_j$, 
with RKHS-norm 
\[
\Big\|\sum_{j\le J} h_j\psi_j\Big\|^2_{\HHH^{s, J}} = \frac1{s^2}\sum_{j \le J} {j^{1+2\alpha}h_j^2}.
\]
The function $f_0$ admits a series expansion $f_0 = \sum f_j\psi_j$. For $J_0 \le J$, 
consider the function $h_0 =  \sum_{j \le J_0} f_j\psi_j$ in the RKHS. It holds that 
\[
\|f_0 - h_0\|^2_2 = \sum_{j > J_0} f^2_k \le J_0^{-2\beta} \|f_0\|^2_\beta.
\]
Hence for $J_0 = \lfloor(\eps/\|f_0\|_\beta)^{-1/\beta}\rfloor$, we have that $\|f_0 - h_0\|_2 \le \eps$. 
The condition on $J$ ensures that $h_0$ is an element of the RKHS, and 
\[
\|h_0\|^2_{\HHH^{s, J}} = \frac1{s^2}\sum_{j \le J_0} {j^{1+2\alpha-2\beta}j^{2\beta}f_j^2}
\le \frac{\|f_0\|^2_\beta}{s^2}J_0^{(1+2\alpha-2\beta) \vee 0}.
\]
It follows that 
\begin{equation}\label{eq:upperboundfortheinfimumoveranepsRKHSballaroundthetrueparameter}
\inf_{h \in \HHH^{s, J} \atop \|h-f_0\|\le \eps}\|h\|^2_{\HHH^{s,J}} \le 
\frac{\|f_0\|^2_\beta}{s^2}J^{{(1+2\alpha-2\beta) \vee 0}}
\end{equation}
Combining this with the preceding lemma and Lemma 5.3 of \cite{rkhs2008}  completes the proof.
\end{prf}

\subsection{Sieves, remaining mass and entropy}

Let the sequence $\eps_n \to 0$ and $\beta > 0$ be given. We consider sieves of growing dimension of the form
\begin{equation}\label{eq: s}
\FF_n = \Big\{h=\sum_{j\le J_n} h_j\psi_j\Big\},
\end{equation}
where 
\begin{equation}\label{eq: j}
J_n = K_1 \eps_n^{-1/\beta}\log 1/\eps_n
\end{equation}
for a constant $K_1 > 0$ specified below. 

By assumption \eqref{eq: p} we have 
\[
\Pi(f \not\in \FF_n ) = \Pi(J > J_n) \lle e^{-C'K_1\eps_n^{-1/\beta} \log 1/\eps_n}. 
\]
This implies that statement \eqref{eq: rm} of Theorem \ref{thm: main} holds if $K_1$ is chosen 
large enough. 

As for the entropy condition \eqref{eq: en}, we note that if the 
function $f_0$ admits the series expansion $f_0 = \sum_{j} f_{0,j}\psi_j$, 
then a function $f \in\FF_n$ which satisfies $\|f-f_0\|_2\le \eps$  is of the form $f = \sum_{j \le J_n} f_j\psi_j$,
and 
$\sum_{j \le J_n} (f_j-f_{0,j})^2 \le \eps^2$.
Hence, the covering number in \eqref{eq: en}  is bounded
by the $a\eps$-covering number of a ball of radius $\eps$ in $\RR^{J_n}$, which is 
bounded by $(3/a)^{J_n}$ (see, for instance, \cite{pollard1990}). 
In view of the choice \eqref{eq: j} of $J_n$ it follows that \eqref{eq: en} holds.

\section{Proof of Theorem \ref{thm: main2}}
\label{sec: proof2}

Under the conditions of Theorem \ref{thm: main2} we can 
replace the result of Theorem \ref{thm: pm1} by the following, 
which implies that \eqref{eq: pm2} holds.

\begin{thm}
Let the prior $\Pi$ on $f$ be defined according to \eqref{eq: p1}--\eqref{eq: p3},  
with $\alpha > 0$ and $p$ and $g$ satisfying \eqref{eq: pp} and \eqref{eq: g},
and let $f_0 \in H^\beta[0,1]$ for $0 < \beta \le \alpha + 1/2$.
Then, for a constant $C > 0$, it holds that  
\[
-\log \Pi(f: \|f-f_0\|_2 \le \eps) \le C \eps^{-1/\beta},
\]
for all $\eps> 0$ small enough.
\end{thm}

\begin{prf}
Instead of using Lemma \ref{lem: cb} we simply note that for $s > 0$ and $J \in \NN$, 
and $Z_1, Z_2, \ldots$ independent and standard normal,
\[
-\log\PP\Big(\Big\|s \sum_{j=1}^J j^{-1/2-\alpha} Z_j \psi_j \Big\|_2 \le \eps\Big)
\le -\log \PP\Big(\Big\| \sum_{j=1}^\infty j^{-1/2-\alpha} Z_j \psi_j \Big\|_2 \le \eps/s\Big).
\]
By Lemma 4.2 of \cite{waaij2016} the right-hand side is bounded by  a constant times $(\eps/s)^{-1/\alpha}$. 
%Combined with  (the proof of) Lemma \ref{lem: sb} 
{ Using  \eqref{eq:upperboundfortheinfimumoveranepsRKHSballaroundthetrueparameter} and Lemma 5.3 of \cite{rkhs2008},  
we see that for $J \ge (\eps/\|f_0\|_\beta)^{-1/\beta}$,} 
\[
-\log \PP\Big(\Big\|s \sum_{j=1}^J j^{-1/2-\alpha} Z_j \psi_j  - f_0\Big\|_2 \le 2\eps\Big) \lle
\Big(\frac s\eps\Big)^{1/\alpha} + 
\frac{1 \vee \eps^{{-(1+2\alpha-2\beta)/\beta}}}{s^2}.
\]
For $\beta \le \alpha + 1/2$ the two terms on the right are balanced for $s$ of the order 
$\eps^{(\beta-\alpha)/\beta}$, 
in which case the right-hand side of bounded by a constant times $\eps^{-1/\beta}$. 
In view of assumption \eqref{eq: pp} and \eqref{eq: g}, it follows by conditioning that 
\[
\Pi(f: \|f-f_0\|_2 \le 2\eps) \ge \exp\big(-c_1\eps^{-1/\beta}\big) p\big(\big\lfloor c_2\eps^{-1/\beta}\big\rfloor\big) 
\int_{\eps^{(2\beta-2\alpha)/\beta}}^{{c_3}\eps^{(2\beta-2\alpha)/\beta}}g(\eta)\,d\eta.  
\]
The assumptions on $p$ and $g$ ensure 
that this is bounded from below by a constant times $\exp(-C\eps^{-1/\beta})${ for some constants \(C,c_1,c_2,c_3>1\).}
\end{prf}

To complete the proof of Theorem \ref{thm: main2} we note that in this case we can use 
the same sieves $\FF_n$ as defined in \eqref{eq: s}, but with a different choice 
for the dimension $J_n$, namely $J_n =\big\lceil K_1\eps_n^{-1/\beta}\big\rceil$, for some $K_1 > 0$. 
The tail condition in \eqref{eq: pp} then ensures that \eqref{eq: rm2} holds if $K_1$
is chosen large enough. The entropy bound \eqref{eq: en2} is obtained by the same
argument as before, but now using the new choice of $J_n$.

\bibliographystyle{harry}
\bibliography{bib}

\end{document}